\documentclass[12pt]{article}
\usepackage{amssymb,latexsym}
\renewcommand{\thefootnote}{}

\begin{document}
\title{On homeomorphisms and quasi-isometries of the real line}
\author{Parameswaran Sankaran\\
Institute of Mathematical Sciences,\\
CIT Campus, Taramani, Chennai 600 113\\ 
Email: {\tt sankaran@imsc.res.in}}
\date{}
\maketitle

\footnote{\noindent A.M.S. Subject Classification (2000):- 
20F65, 20F28, 20F67\\
Key words and phrases: pl-homeomorphisms,   
quasi-isometry, Thompson's group, free groups.}

\thispagestyle{empty}

\def\theequation {\arabic{section}.\arabic{equation}}
\renewcommand{\thefootnote}{}

\newcommand{\codim}{\mbox{{\rm codim}$\,$}}
\newcommand{\stab}{\mbox{{\rm stab}$\,$}}
\newcommand{\lr}{\mbox{$\longrightarrow$}}

\newcommand{\ch}{{\cal H}}
\newcommand{\cf}{{\cal F}}
\newcommand{\cd}{{\cal D}}

\newcommand{\blr}{\Big \longrightarrow}
\newcommand{\da}{\Big \downarrow}
\newcommand{\ua}{\Big \uparrow}
\newcommand{\hra}{\mbox{\LARGE{$\hookrightarrow$}}}
\newcommand{\rt}{\mbox{\Large{$\rightarrowtail$}}}
\newcommand{\dua}{\begin{array}[t]{c}
\Big\uparrow \\ [-4mm]
\scriptscriptstyle \wedge \end{array}}

\newcommand{\be}{\begin{equation}}
\newcommand{\ee}{\end{equation}}
\newcommand{\ov}{\overline}
\newtheorem{guess}{Theorem}[section]
\newcommand{\bth}{\begin{guess}$\!\!\!${\bf .}~}
\newcommand{\eeth}{\end{guess}}
\renewcommand{\bar}{\overline}
\newtheorem{propo}[guess]{Proposition}
\newcommand{\bpropo}{\begin{propo}$\!\!\!${\bf .}~}
\newcommand{\epropo}{\end{propo}}

\newtheorem{lema}[guess]{Lemma}
\newcommand{\blem}{\begin{lema}$\!\!\!${\bf .}~}
\newcommand{\elem}{\end{lema}}

\newtheorem{defe}[guess]{Definition}
\newcommand{\bdefe}{\begin{defe}$\!\!\!${\bf .}~}
\newcommand{\edefe}{\end{defe}}

\newtheorem{coro}[guess]{Corollary}
\newcommand{\bcor}{\begin{coro}$\!\!\!${\bf .}~}
\newcommand{\ecor}{\end{coro}}

\newtheorem{rema}[guess]{Remark}
\newcommand{\brem}{\begin{rema}$\!\!\!${\bf .}~\rm}
\newcommand{\erem}{\end{rema}}

\newtheorem{exam}[guess]{Example}
\newcommand{\beg}{\begin{exam}$\!\!\!${\bf .}~\rm}
\newcommand{\eeg}{\end{exam}}

\newcommand{\ctext}[1]{\makebox(0,0){#1}}
\setlength{\unitlength}{0.1mm}
\newcommand{\cl}{{\cal L}}
\newcommand{\cp}{{\cal P}}
\newcommand{\ci}{{\cal I}}
\newcommand{\bz}{\mathbb{Z}}
\newcommand{\cs}{{\cal s}}  
\newcommand{\cv}{{\cal V}}
\newcommand{\ce}{{\cal E}}
\newcommand{\ck}{{\cal K}}
\newcommand{\cR}{{\cal R}}
\newcommand{\bq}{\mathbb{Q}}
\newcommand{\bt}{\mathbb{T}}
\newcommand{\bh}{\mathbb{H}}
\newcommand{\br}{\mathbb{R}}
\newcommand{\wt}{\widetilde}
\newcommand{\im}{{\rm Im}\,}
\newcommand{\bc}{\mathbb{C}}
\newcommand{\bp}{\mathbb{P}}
\newcommand{\bn}{\mathbb{n}}
\newcommand{\ds}{\displaystyle}
\newcommand{\tor}{{\rm Tor}\,}
\newcommand{\bs}{\mathbb{S}}
\def\ns{\mathop{\lr}}
\def\nssup{\mathop{\lr\,sup}}
\def\nsinf{\mathop{\lr\,inf}}
\renewcommand{\phi}{\varphi}
\newcommand{\co}{{\cal O}}
\newcommand{\pl}{PL_\delta(\br)}
\newcommand{\plk}{PL_\kappa(\br)} 
\newcommand{\plo}{PL_\delta^+(\br)}
\newcommand{\pls}{PL(\bs^1)}
\newcommand{\plsr}{\wt{PL}(\bs^1)}
\newcommand{\qr}{QI(\br)}
\newcommand{\qro}{QI^+(\br)}
\newcommand{\qz}{QI(\bz)}
\newcommand{\supp}{\mbox{\rm Supp}}
\noindent
{\bf Abstract:} We show that the group of piecewise-linear homeomorphisms 
of $\br$ having bounded slopes surjects onto the group $QI(\br)$
of all quasi-isometries of $\br$. We prove that the 
following groups can be imbedded in $QI(\br)$: the  
group of compactly supported piecewise-linear homeomorphisms of 
$\br$, the Richard Thompson group $F$,  
and the free group of continuous rank.   
    
\section{Introduction}
We begin by recalling the notion of quasi-isometry. 
Let $f:X\lr X'$ be a map (which is not assumed to be continuous)  
between metric spaces. We say that $f$ is a $C$-quasi-
isometric embedding if there exists a $C>1$ such that 
$$C^{-1}d(x,y)-C \leq d'(f(x),f(y)) 
\leq C d(x,y)+C \eqno(*)$$ for all $x,y\in X$. Here $d,d'$ denote the 
metrics on $X,X'$ respectively. 
If, further,  every $x'\in X'$ is 
within distance $C$ from the image of $f$, we say that $f$ is a 
$C$-quasi isometry. If $f$ is a quasi-isometry (for some $C$) 
then there exists a quasi-isometry $f':X'\lr  X$ (for a possibly  
different constant  $C'$) such that 
$f'\circ f $ (resp. $f\circ f'$)  is quasi-isometry equivalent to the 
identity map of $X$  (resp. $X'$). (Two maps $f,g: X\lr X$ are said to be  
quasi isometrically  equivalent if there exists a constant $M$ such 
that $d(f(x),g(x))\leq M$ for all $x\in X$.) Let $[f]$ denote the 
equivalence class of a quasi-isometry $f:X\lr X$. 
The set $QI(X)$ of all equivalence classes of quasi-isometries of $X$ 
is a group under composition: $[f].[g]=[f\circ g]$ for $[f],[g]\in QI(X)$.
If $X'$ is quasi-isometry equivalent to $X$, then $QI(X')$ is 
isomorphic to $QI(X)$. We refer the reader to \cite{bh} for 
basic facts concerning quasi-isometry. For example $t\mapsto [t]$ 
is a quasi-isometry from $\br$ to $\bz$.  

Let $f:\br \lr \br$ be any homeomorphism of $\br$. Denote by 
$B(f)$ the set of break points of  $f$, i.e., points where $f$ 
fails to have 
derivative and by $\Lambda(f)$ the set of slopes of $f$, i.e., 
$\Lambda(f)=\{f'(t)\mid t\in \br\setminus B(f)\}$. 
Note that 
$B(f)\subset\br$ is discrete if $f$ is piecewise differentiable. 

\bdefe\label{pld}
We say that a subset $\Lambda$ of $\br^*$, the set of non-zero real numbers,  
is {\em bounded} 
if there exists an $M>1$ 
such that $M^{-1}<|\lambda|<M$ for all $\lambda\in \Lambda$.  
We say that a homeomorphism $f$ of  $\br$ which is piecewise 
differentiable has {\em bounded slopes} if $\Lambda(f)$ is bounded. 

We denote by $PL_\delta(\br)$ the set of all those   
piecewise-linear homeomorphisms $f$ of  $\br$ such that   
$\Lambda(f)$ is bounded.  It is clear that 
$\pl$ is a subgroup of the group $PL(\br)$ of all piecewise-linear 
homeomorphisms 
of $\br$. 
\edefe 

It is easy to see that 
each $f\in PL_\delta (\br)$ is a quasi-isometry. (See lemma 
\ref{basic} below.)  
One has a 
natural homomorphism $\phi:\pl\lr QI(\br)$ where $\phi(f)=[f]$ for 
all $f\in \pl$. 

\bth \label{main} 
The natural homomorphism $\phi:\pl\lr QI(\br)$, defined as 
$f\mapsto [f]$,  is surjective.
\eeth 

If $f:\br\lr \br$ is a homeomorphism, recall that 
$\supp(f) $, the support of $f,$ is the closure 
of the set $\{x\in \br\mid f(x)\neq x\}$ of all points 
moved by $f$.
Denote by $PL_\kappa(\br)$ the group of all 
piecewise-linear homeomorphisms of $\br$ which 
have compact support. It is obvious 
that $PL_\kappa(\br)\subset\ker(\phi).$

Let $\Gamma$ be a group of homeomorphisms of $\bs^1$. 
Any $f\in \Gamma$ can be lifted to obtain a homeomorphism 
$\wt{f}$ of $\br$ over the covering projection 
$p:\br \lr \bs^1$, $t\mapsto \exp(2\pi\sqrt{-1}t)$. 
The set $\wt{\Gamma}$  of all homeomorphisms of $\br$ which 
are lifts of elements of $\Gamma$ is a subgroup of 
the group $Homeo(\br)$ of all homeomorphisms of $\br$. 
Indeed $\wt{\Gamma}$ is a central extension of 
$\Gamma$ by the infinite cyclic group generated by 
translation by $1$: $x\mapsto x+1$.  Denote by $Diff(\bs^1)$ 
the group of all $C^\infty$ diffeomorphisms of the circle. 
When $\Gamma$ is one of 
the groups $\pls$, $Diff(\bs^1)$, any element of 
$\wt{f}\in \wt{\Gamma}$  has bounded slope and  
is quasi-isometrically equivalent to the identity map of $\br$ 
(since $\wt{f}(x+n)=\wt{f}(x)+n$ for $n\in \bz$). 

Recall that Richard Thompson discovered the group 
$$F=\langle x_0,x_1,\cdots \mid x_ix_jx_i^{-1}=x_{j+1},~i<j\rangle$$ 
and used it in some constructions in logic related to  
word problems. The group $F$ is finitely presentable with 
two generators $x_0,x_1$ and two relations. 
This group and a closely related larger group 
$G$ have since then appeared in several contexts including 
homotopy theory \cite{fh}, 
homological group theory \cite{bg}, Teichm\"uller 
theory \cite{im}, etc. The group $F$ is isomorphic to  
the subgroup of piecewise-linear homeomorphisms of $\br$ which 
are the identity outside the unit interval $I$ such 
that $B(f)$ is contained in dyadic rationals and 
$\Lambda(f)$ is contained in the subgroup of $\br^*$ 
generated by $2$. Although $F$ satisfies no (nontrivial) group 
law, it contains no non-abelian free group. The group $G$ 
is the group of piecewise-linear  
homeomorphisms $f$ of the circle $\bs^1=I/\{0,1\}$ with 
$B(\wt{f})$ contained in dyadic rationals and $\Lambda(\wt{f})$ contained 
in the multiplicative subgroup of $\br^*$ generated by $2$ for some 
lift $\wt{f}$ of $f$. It 
is the first known example of a finitely presented infinite simple 
group. We recommend the beautiful survey article 
\cite{cfp} for further information about Richard Thompson's groups. 

We shall prove the following theorem: 

\bth \label{free} 
The following groups can be imbedded in $QI(\br)$.\\
(i) the groups $\wt{Diff}(\bs^1)$ and $\plsr$, \\ 
(ii)  the group $PL_\kappa(\br)$, \\      
(iii)  the Thompson's group $F$, and, \\  
(iv) the free group of rank $c$, the continuum. 
\eeth 

Our proofs are completely elementary.           
We explain the main idea of the proof of theorem 
\ref{free}. Take for example the group $\plk$. 
The first step is to realise this as a subgroup $\Gamma_1$
of $\plk$ having support in $(0,1)$. 
This is achieved easily by imbedding 
$\br$ in the interval $(0,1)$.  The group  
$\Gamma_1$ can be thought of as a group of 
piecewise-linear homeomorphisms of the circle. Lifting this 
back to $\br$ via the covering projection, 
we obtain now a group 
$\wt{\Gamma}_1$ which no longer has compact 
support. However each element of this 
group is quasi-isometric to $id$. So we conjugate 
this group by a piecewise-linear homeomorphism whose slope 
grows exponentially. The result is that the 
features of each element of $\wt{\Gamma}_1$ get 
magnified resulting in a quasi-isometry not 
representing $1$.  The same trick works for 
$\wt{Diff}(\bs^1)$ as well. 
Parts (iii) and (iv) follow from known 
embeddings of the relevant groups.

%%%%%%%%%%%%%%%%%%%%%%%%%%%%%%%%%%%%%%%%%%%%%%%%  
%%%%%%%%%%%%%     PROOFS     %%%%%%%%%%%%%%%%%%%
%%%%%%%%%%%%%%%%%%%%%%%%%%%%%%%%%%%%%%%%%%%%%%%%

\section{Proof of Theorem \ref{main}}  
We first establish the following basic observation.

\blem \label{basic} 
Let $f$ be a piecewise differentiable homeomorphism of $\br$ 
with $\Lambda(f)\subset \br^*$ bounded.   
Then $f$ is a quasi-isometry. 
\elem  
\noindent 
{\bf Proof:} Replacing $f$ by $-f$ if necessary, one may assume 
without loss of generality that $f$ is monotone increasing. 

Suppose that $\Lambda(f)\subset (1/M, M)$. 
If 
$f$ is differentiable everywhere, then it is an 
$M$-quasi-isometry. 

Suppose that $B(f)\neq \emptyset$.  Let $a\in \br$.   
For any $b>a$, let $a_1<\cdots <a_k$ be the points 
of $(a,b)$ where $f$ is non-differentiable.
%and let $\lambda_i$ 
%be the slope of $f$ in $(a_i,a_{i+1})$ where $a_0:=a, a_{k+1}:=b$. 
Then, applying the mean value theorem, 
$f(b)-f(a)=\sum_{0\leq i\leq k}(f(a_{i+1})-f(a_i))
=\sum_{0\leq i\leq k} f'(c_i)(a_{i+1}-a_i)$ for some $c_i\in (a_i,a_{i+1})$. 
Since $\Lambda(f)\subset (1/M,M)$, it follows that 
$M^{-1}(b-a)<f(b)-f(a)<M(b-a)$.  Since 
$a,b\in \br$ are arbitrary,  we conclude that 
$f$ is a quasi-isometry. \hfill $\Box$ 

One has a well-defined map $\phi:\pl \lr QI(\br)$ which 
is a homomorphism. We now prove that $\phi$ is surjective. 

\blem \label{monot} 
Let $f:\br \lr \br$ be a $C$-quasi-isometry that preserves the ends 
of $\br$.  Let $x\in \br$. Then (i)  there exists $y$ such that 
$y-x\leq 4C^2$ is positive integer and $f(y)>f(x)$; (ii) there 
exists $v$ such that $x-v\leq 4C^2$ is a positive 
integer and $f(x)>f(v)$.  
\elem 
\noindent
{\bf Proof:} If $f(x+1)>f(x)$, then $y=x+1$ meets our requirements. 

Assume that $f(x)>f(y)$ for all $y$ such that 
$x+1\leq y<4C^2+x$. Let $z\geq x+2$ be the smallest real number such 
that $z-x$ is a positive integer and 
$f(z)>f(x)\geq f(z-1)$. Such a $z$ exists since $f(t)\rightarrow +\infty$  
as $t\rightarrow +\infty$. By our assumption $z-x\geq 4C^2+1.$  
Set $u=z-1$.  Then the inequality ($*$) implies  
$f(u)<f(x) +C -C^{-1}(u-x)\leq f(x)-3C$ and $f(z)-f(u)<C(z-u)+C=2C$. 
Hence $f(z)<f(u)+2C<f(x)-C$, i.e., $f(z)-f(x)<-C$. This contradicts 
our hypothesis that $f(z)>f(x)$, completing the proof of part (i).
Proof of part (ii) is similar. \hfill $\Box$

\noindent 
{\bf Proof of Theorem \ref{main}:}  
Since the subgroup  $QI^+(\br)\subset QI(\br)$ 
that preserves the ends $\{+\infty,-\infty\}$ of $\br$ is of index $2$ 
and since $\pl$ contains 
elements which are orientation reversing, it suffices 
to show that $\qro$ is contained in the image of $\phi$ where  
$\qro\subset \qr$ is the index $2$ subgroup whose elements preserve 
the ends of $\br$.

Let $f:\br\lr \br$ be a $C$-quasi isometry, with $C>1$,  
which preserves 
the ends of $\br$. We assume, as we may, that $C$ is a positive integer.   

Set $x_0=0$. We define $x_k\in \bz$ for any  
integer $k$ as follows: Let $k\geq 1$. Having defined 
$x_{k-1}$ inductively, choose $x_k>x_{k-1}$ to be the 
smallest integer such that $f(x_k)>f(x_{k-1})$.   
For any negative integer $k$, we define $x_k$ analogously 
(by downward induction) as the greatest integer such that $x_k<x_{k+1}$ and 
$f(x_k)<f(x_{k+1})$.  

Set $y_k:=x_{C^3k}$, and let 
$B:=\{y_k| ~k\in\bz\}\subset \bz$.  By lemma \ref{monot}, 
we see that $B$ is a discrete subset of $\br$ which is  
$4C^5$-dense in $\br$.  Note that for any $k\in \bz$, 
$y_k-y_{k-1}\geq C^3$. 

Since $f(y_k)>f(y_{k-1})$ for all $k\in \bz$, there exists a unique 
piecewise-linear homeomorphism  $g:\br\lr\br$ such that $g(y_k)=f(y_k)$  
and is {\it linear} on the interval $[y_{k-1}, y_k]$ for every $k\in \bz$.  
We claim that $g$ has bounded slopes. 
Since $g$ is linear on each of the intervals $[y_{k-1},y_k]$, 
we need only bound $\frac{g(y_{k})-g(y_{k-1})}{y_k-y_{k-1}}$.  Indeed, 

$$\frac{g(y_k)-g(y_{k-1})}{y_k-y_{k-1}}=\frac{f(y_k)-f(y_{k-1})}{y_k-y_{k-1}}
<C+\frac{C}{y_k-y_{k-1}} \leq C+C^{-2} $$  
as $y_k-y_{k-1}\geq C^3$. 
Similarly,\\  
 $$\frac{g(y_k)-g(y_{k-1})}{y_k-y_{k-1}}
>C^{-1}-C^{-2}.$$
 
It follows that $\Lambda(g)\subset [C^{-1}-C^{-2}, C+C^{-2}]$ and  
$g\in \pl$.  
 
Since $f$ and $g$ agree on the quasi-dense set $B$, 
we see that $[f]=[g]$. This completes the proof. \hfill $\Box$ 

\brem 
\noindent 
(i)  By setting $g(y_k)$ equal to a rational number sufficiently 
close to $f(y_k)$ in the above proof, we see that since 
$y_k\in \bz$, the element $g\in \pl$ has rational slopes.     
Consequently it follows that $\phi$ restricted to the subgroup 
$PL_\bq^{\bq^*}(\br)$ of $\pl$ 
consisting of those $g\in \pl$ having slopes in $\bq^*$ and $B(g)$ 
contained in $\bq$ is surjective.\\ 

(ii) The kernel of $\phi$ contains the group of 
all piecewise-linear homeomorphisms which have slope $1$ 
outside a compact interval.  This latter group equals to 
the derived group $PLF'(\br)$ where $PLF(\br)$ denotes 
the subgroup of $\pl$ consisting of homeomorphisms $f$ for 
which $B(f)$ is finite.  Also $\plk=PLF''(\br)$. See \cite{bs}.   
\erem 

\section{Proof of Theorem \ref{free}}

Let $h_1:\br\lr (0,1)$ be the homeomorphism defined by 
$h_1(-x)=1-h_1(x)$ for every $x\in\br,~h_1(n)=1-1/(n+2)$ for each 
integer $n\geq 0$ and is linear on 
each interval $[n,n+1]$ for $n\in \bz$.  If $f$ is any 
compactly supported (piecewise-linear) homeomorphism of $\br$ 
then $h_1\circ f\circ h_1^{-1}$ 
is a compactly supported (piecewise-linear) homeomorphism of $(0,1)$.   
Since $\bs^1=I/\{0,1\}$, we also get an embedding
$\bar{\eta}:\plk\lr PL(\bs^1)$ where $\bar{\eta}(f)$
is defined to be the extension of $h_1\circ f\circ h_1^{-1}$ 
to $\bs^1$. We define 
$\eta:\plk\lr \pl$ as the imbedding $f\mapsto \eta(f)$ where 
$\eta(f)(n)=n$ for $n\in \bz$ and 
$\eta(f)(x)=n+h_1fh_1^{-1}(x-n)$ for $n< x<n+1$.  

Let $h_0:\br\lr \br$ be the piecewise-linear homeomorphism defined as 
follows: $h_0(-x)=-h_0(x)~\forall x\in\br, h_0(x)=x$ for
$0\leq x\leq 1$
and maps the interval $[n,n+1]$ onto $[2^{n-1},2^n]$
linearly for each positive integer $n$. 

Suppose $f:\bs^1\lr\bs^1$ is an orientation preserving 
piecewise-linear homeomorphism or a diffeomorphism. Let $\wt{f}:\br\lr\br$ be 
any lift of $f$ so that $p\circ\wt{f}=f\circ p$, where 
$p:\br \lr \bs^1$ is the covering projection $t\mapsto \exp(2\pi\sqrt{-1}t)$.  
Then $[\wt{f}]=1$ in $\qr$. (Indeed one has 
$\wt{f}(x+n)=n+\wt{f}(x)$ for all $x\in \br$ and $n\in \bz$ and so   
$|\wt{f}-id|\leq |\wt{f}(0)|+1$.) 

Let $\Gamma$  be one of the groups $\pls$  or $Diff(\bs^1)$ and 
let $\wt{\Gamma}$ be the group of homeomorphisms of $\br$ 
which are lifts of elements of $\Gamma$ with respect to 
the covering projection $p$.  For $\wt{f}\in \wt{\Gamma}$ 
set $f_0:=h_0\wt{f}h_0^{-1}$.  
Clearly, $\wt{f}\mapsto f_0$ is a monomorphism of groups 
$\wt{\Gamma}\lr Homeo(\br)$.   
We claim that for any $\wt{f}\in\wt{\Gamma}$, 
$f_0$ is a quasi-isometry.  To see this, we assume 
without loss of generality that $\wt{f}$ is orientation 
preserving.  It is clear that $f_0$ is 
differentiable outside a discrete subset 
of $\br$. 
We claim that $f_0$ has bounded slopes. Since 
$f_0$ has continuous 
derivatives on each interval on which $f_0$ has derivatives, 
it suffices to show that the set $\{f_0'(t)\}$ as $t$ varies in 
$\br \setminus B$ 
is bounded, where $B$ is any discrete set which contains 
$B(f_0)$.  We set $B:=B(h_0)\cup h_0B(\wt{f})\cup h_0\wt{f}^{-1}B(h_0)$. 

Let $0<m<M$ be such that $m<\wt{f}'(x)<M$ for $x\in \br$.
Let $t\in \br\setminus B$ and set $s=h_0^{-1}(t), u=\wt{f}(s)$ so that  
$h_0^{-1}, \wt{f}, h_0$ are  differentiable at  $t, s, u$ respectively.  
Consequently $f_0$ is differentiable at $t$.  

Since $|u-s|=|\wt{f}(s)-s|<q$ where $q:=[|\wt{f}(0)|]+2$,  
we see that $2^{-q}<h_0'(u)/h'_0(s)<2^q$.  
Using the chain rule, it follows  
that $f_0'(t)=h_0'(u)\wt{f}'(s) (h_0^{-1})'(t)=\wt{f}'(s)h_0'(u)/h_0'(s)$
lies in the interval $(2^{-q}m,2^qM)$. 
It follows from lemma \ref{basic} that  
$f_0$ is a quasi-isometry. 

It is clear that the  map $\psi: \wt{\Gamma} \lr \qr$ defined 
as $ \wt{f}\mapsto [f_0]$ is a homomorphism.   

We are now ready to prove theorem \ref{free}.\\

\noindent 
{\bf Proof of theorem \ref{free}:} We use the above notations 
throughout the proof. 

\noindent
(i) We prove that $\psi:\wt{\Gamma}\lr \qr$ is a monomorphism
where $\Gamma=\pls$ or $Diff(\bs^1)$.
Suppose that $\wt{f}\in \wt{\Gamma}$,
$\wt{f}\neq id$.     
We shall show that that $|f_0-id|$ is unbounded.  
Choose $x$ in the interval $[0,1)$ such that $\wt{f}(x)\neq x$. 
Set $k=[\wt{f}(x)]$ so that $\wt{f}(x)=k+y$, $0\leq y<1$.  Replacing $\wt{f}$ 
by its inverse if necessary, we  assume without loss of generality that 
$x<\wt{f}(x)$. This implies that $k\geq 0$ with equality only if $y>x$.
For any positive integer $n$, we have 
$f_0(2^n+2^nx)=h_0\wt{f}h_0^{-1}(2^n+2^nx)
=h_0\wt{f}(n+1+x)=h_0(n+1+\wt{f}(x))=h_0(n+1+k+y)
=2^{n+k}+2^{n+k}y$.\\ 
If $k=0$, then $y>x$ and so $f_0(2^n+2^nx)-(2^n+2^nx)=2^n(y-x)$. Thus 
$|f_0-id|$ is unbounded. \\  
If $k>0$, then $f_0(2^n+2^nx)-(2^n+2^nx)=2^{n+k}+2^{n+k}y-2^n-2^nx 
\geq 2^{n+1}-2^n-2^nx=2^n(1-x)$.   
As $0\leq x<1$, again it follows that $|f_0-id|$ is unbounded.

\noindent
(ii) As observed earlier, 
$\eta:\plk\lr \pl$  is a monomorphism. It is evident that 
the image of $\eta$ is contained in $\plsr$. Since
$\psi:\plsr \lr \qr$ is a monomorphism by (i), 
assertion (ii) follows.

\noindent 
(iii) Now statement (iii) follows from (ii) above and the fact that
Thompson's group $F$ is isomorphic to the subgroup of $\plk$ 
of all piecewise-linear homeomorphisms which have support in  
$[0,1]$ having break points contained in the
set of dyadic rationals in $[0,1]$ and slopes contained 
in the multiplicative subgroup of $\br^*$ generated by $2$.

\noindent
(iv) To prove (iv), recall that Grabowski \cite{gra} has shown that 
the free group of rank $c$ the continuum 
embeds in the group of compactly supported $C^k$  
diffeomorphisms ($1\leq k\leq \infty$) 
of any positive dimensional manifold. 
In particular this is true of $Diff(\bs^1)$. It follows  
easily that $\wt{Diff}(\bs^1)$ also contains a 
free group of rank the continuum. By part (i), 
this completes the proof. \hfill $\Box$ 

\blem \label{torsion}
The group $\qro$ is torsion-free.  
\elem 
\noindent
{\bf Proof:} Let $f\in \pl$ be such that $[f]\neq 1\in \qro$.
Thus $f-id$ is unbounded. Choose a sequence $(a_n)$ 
of real numbers such that  
$a_n\rightarrow +\infty$ as $n\rightarrow +\infty$  and 
$|f(a_n)-a_n|\rightarrow +\infty$.  
Let $k>1$ be any integer.  Suppose that  
$f(a_n)>a_n$.   Since $f$ is order preseving, 
for each $n$ we have $a_n<f(a_n)<\cdots <f^k(a_n).$  
In particular $f^k(a_n)-a_n>f(a_n)-a_n.$ 
Similarly,  $a_n-f^k(a_n)>a_n-f(a_n)$ in case  
$a_n>f(a_n)$. Therefore $|f^k(a_n)-a_n|>|f(a_n)-a_n|~\forall n$ 
and hence $f^k-id$ is unbounded. Hence 
$[f^k]\neq 1$ in $\qro$ for $k>1$. 
\hfill $\Box$ 
  
\brem  Thompson's group $G$ does not imbed in $\qr$
since it has an element of order $3$ whereas it follows from 
Lemma \ref{torsion} that 
all torsion elements in $\qr$ are of order $2$.  
\erem 

{\bf Acknowledgements:} 
Part of this work was done while the author was visiting
the University of Calgary, Alberta, Canada, during the
Spring and Summer of 2003. 
It is a pleasure to thank Professors K.Varadarajan and P.Zvengrowski
for their invitation and hospitality as well as
financial support through their NSERC grants making
this visit possible.

%%%%%%%%%%%%%%%%%%%%%%%%%%%%%%%%%%%%%%%%%%%%%%%%%%%%%%%%%%%%%%%%  

\end{document}